\documentclass[reqno,12pt]{amsart}
\usepackage{amssymb,amsmath,amsthm,mathtools,mathrsfs,calc,verbatim}
\usepackage{enumerate}
\usepackage{bbm}
\usepackage[top=1.5in, bottom=2in, left=1.5in, right=1.5in]{geometry}
\usepackage{filecontents}

\def\Z{\mathbb{Z}}

\def\H{\mathbb{H}}
\def\N{\mathbb{N}}
\def\C{\mathbb{C}}

\DeclareMathOperator{\Li}{Li}

\newcommand{\pfrac}[2]{\left(\frac{#1}{#2}\right)}

\renewcommand{\bar}[1]{\overline{#1}}

\newtheorem{theorem}{Theorem}
\numberwithin{theorem}{section}

\theoremstyle{remark}
\newtheorem*{remark}{Remark}

\numberwithin{equation}{section}
\theoremstyle{definition}
\newtheorem*{definition}{Definition}

\makeatletter
\renewcommand{\pod}[1]{\mathchoice
  {\allowbreak \if@display \mkern 18mu\else \mkern 8mu\fi (#1)}
  {\allowbreak \if@display \mkern 18mu\else \mkern 8mu\fi (#1)}
  {\mkern4mu(#1)}
  {\mkern4mu(#1)}
}

\begin{document}

\title{Images of Maass-Poincar\'e series in the lower half-plane}

\thanks{The research of the second author is supported by the Alfried Krupp Prize for Young University Teachers of the Krupp foundation and the research leading to these results receives funding from the European Research Council under the European Union's Seventh Framework Programme (FP/2007-2013) / ERC Grant agreement n. 335220 - AQSER}

\author{Nickolas Andersen, Kathrin Bringmann, and Larry Rolen}

\address{Mathematics Department, UCLA, Box 951555, Los Angeles, CA 90095}
\email{nandersen@math.ucla.edu}

\address{Mathematical Institute, University of Cologne, Weyertal 86-90, 50931 Cologne, Germany}
\email{kbringma@math.uni-koeln.de}

\address{Hamilton Mathematics Institute \& School of Mathematics, Trinity College, Dublin 2, Ireland}
\email{lrolen@maths.tcd.ie}

\begin{abstract}
In this note we extend integral weight harmonic Maass forms to functions defined on the upper and lower half-planes using the method of Poincar\'e series.
This relates to Rademacher's ``expansion of zero'' principle, which was recently employed by Rhoades to link mock theta functions and partial theta functions.
\end{abstract}

\maketitle

\date{\today}


\section{Introduction and statement of results}
In \cite{Rh}, Rhoades found a method to uniformly describe partial theta functions and mock theta functions as manifestations of a single function.
He showed that Ramanujan's mock theta function $f(q)$ (defined below), with $q:=e^{2\pi i\tau}$ and $\tau$ in the upper half-plane $\H$, in some sense ``leaks'' through the real line to a partial theta function $\psi(q^{-1})$ (given below) on the lower half-plane $-\H$. 
His construction follows the ``expansion of zero'' principle of Rademacher (see \cite{Knopp}, \cite[Chapter IX]{Leh-Disc}, and  \cite{Rad-PNAS}).
Rademacher showed, using his exact formula for the partition function, that the partition generating function can be extended to the lower half-plane, and he later proved \cite{Rad-Top} that this extension is identically zero in the lower half plane. We note that there are other relations of partial theta functions and mock theta functions. For example one, which is due to Zagier and Zwegers, passes through asymptotic expansions (see for example \cite{LZ}).

Let us now say a few more words concerning mock theta functions. Originally introduced by Ramanujan in his last letter to Hardy, mock theta functions have since found applications in many areas of mathematics.
We now understand that they fit into the larger framework of harmonic Maass forms, as shown by Zwegers \cite{Zw-thesis} (see also \cite{Ono-visions,Zag-mock}).
That is, the mock theta functions are examples of mock modular forms, which are the holomorphic parts of harmonic Maass forms (see Section \ref{sec:preliminaries} for definitions).
Thus, it is natural to ask whether Rhoades' construction applies to the non-holomorphic completion of $f(q)$ and, if so, what is the image of that function in the lower half-plane?
One hope, which has not yet been realized, is that this might shed some light on the problem of finding a completion of the partial theta functions to non-holomorphic modular forms.
General {\it partial theta functions} have the shape
\[
\sum_{n\ge 0}\psi(n) n^\nu q^{n^2}
\]
with $\psi$ a primitive Dirichlet character and $\nu\in\Z$ such that $\psi(-1)=(-1)^{\nu+1}$. In particular these functions are not modular forms.

We begin by more closely recalling Rhoades' results.
The partial theta function which Rhoades studied is given by
$$
\psi(q):=\sum_{n\geq1}\left(\frac{-12}{n}\right)q^{\frac{n^2-1}{24}},
$$
and the associated mock theta function is Ramanujan's third order function
$$
f(q):=\sum_{n\geq0}\frac{q^{n^2}}{(1+q)^2\cdot\ \dots\ \cdot(1+q^n)^2}.
$$
Now set
$$
\alpha_c(s):=\sum_{m\geq0}\left(\frac{\pi}{12c}\right)^{2m+\frac12}\frac{1}{\Gamma\left(m+\frac32\right)}\frac{1}{s^{m+1}}
$$
and (with $\zeta_a^b:=e^{\frac{2\pi ib}{a}}$)
$$
\Phi_{c,d}(\tau):=\frac{1}{2\pi i}\int_{|s|=r}\frac{\alpha_c(s)e^{23s}}{1-\zeta_{2c}^d \, q \, e^{24s}}ds,
$$
where $r$ is taken sufficiently small such that $|\operatorname{Log}\left(\zeta_{2c}^dq\right)|\gg r$ and such that the integral converges. Moreover let $\omega_{h,c}$ be the multiplier of the Dedekind $\eta$-function (which can be given explicitly in terms of Dedekind sums, see \cite{Rad-Top}). Then define the function
\begin{align*}
F(\tau)&:=1+\pi\sum_{c\geq1}\frac{(-1)^{\left\lfloor \frac{c+1}{2}\right\rfloor}}{c}\\
&\times\sum_{d\pmod{2c}^\ast}\omega_{-d,2c}\exp\left(2\pi i\left(-\frac{d}{8}\left(1+(-1)^c\right)+\frac{d}{2c}+\tau\right)\right)\Phi_{c,d}(\tau),
\end{align*}
where $d\pmod{2c}^*$ indicates that the sum ranges over those $d$ modulo ${2c}$ with $\gcd(d,2c)=1$.
This function converges in both the upper and lower half-planes, i.e., for $\tau\in\mathbb H\cup (-\mathbb H)$. Moreover, Rhoades' main result states that
$$
F(\tau)=\begin{cases} f(q)\qquad&\text{if }\tau\in\mathbb{H},\\
2\psi\left(q^{-1}\right)&\text{if }\tau\in-\mathbb H.
\end{cases}
$$

As discussed above, we describe a similar phenomenon for both the holomorphic and non-holomorphic parts of Maass--Poincar\'e series.
To state our results, we first require some notation.
Throughout, let $k\in 2\Z$, and let $M_k^!(\Gamma_0(N))$ denote the space of weakly holomorphic modular forms of weight $k$ on $\Gamma_0(N)$.
Let $S_k^!(\Gamma_0(N))$ denote the subspace of $M_k^!(\Gamma_0(N))$ consisting of forms whose Fourier expansion at $i\infty$ has constant term equal to zero.
For $f(\tau)=:\sum_n c_f(n) q^n \in S_k^!(\Gamma_0(N))$, we define the (holomorphic) {\it Eichler integral}
\begin{equation*}
	\mathcal{E}_f(\tau) := \sum_{n\in\Z\backslash\{0\}} \frac{c_f(n)}{n^{k-1}} q^n
\end{equation*}
and the {\it non-holomorphic Eichler integral} ($\tau=u+iv$)
\begin{equation*}
f^\ast(\tau) := -(4\pi)^{1-k}\sum_{n\in\Z\setminus\{0\}} \frac{\bar{c_f(n)}}{n^{k-1}} \, \Gamma(k-1,4\pi nv) q^{-n}.
\end{equation*}
Here $\Gamma(s, y)$ denotes the incomplete gamma function defined in \eqref{incomplete}. Note that with
\[
	D^{k-1} := \left(q\frac{d}{dq}\right)^{k-1} \qquad \text{and} \qquad  \xi_k:=2iv^{k} \frac{\overline{\,\partial\,}}{\partial\overline{\tau}},
\]
we have
\[
	D^{k-1}(\mathcal E_f) = f \qquad \text{and} \qquad  \xi_{2-k}(f^\ast)=f.
\]


For even $k>2$ and $m\in \Z$, let $P_{k,m}$ denote the holomorphic Poincar\'e series (defined in Section~\ref{sec:preliminaries} below). If $m<0$, these functions are weakly holomorphic forms, while for $m>0$, they are cusp forms.
For $k\in-2\N$ and $m>0$, let $F_{k,-m}$ be the Maass-Poincar\'e series  of Section~\ref{sec:preliminaries}. They are harmonic Maass forms with exponential growth in their holomorphic part.

\begin{theorem} \label{thm:integer-weight}
Let $k\in -2\N$ and $m\in\N$. Then the function $H_{k,m}:= H_{k,m}^+ + H_{k,m}^-$ \textup(defined in \eqref{eq:H-k-m-+} and \eqref{eq:h-k-m--} below\textup) converges for all $\tau\in \H\cup (-\H)$. Furthermore, if $\tau\in \H$ we have
\begin{equation*}
	H_{k,m}(\tau) = F_{k,-m}(\tau),
\end{equation*}
and if $\tau\in-\H$ we have 
\begin{equation*}
	H_{k,m}(\tau) = m^{1-k}\left(\mathcal E_{P_{2-k,m}} (-\tau)- \frac{(4\pi)^{1-k}}{(-k)!}P_{2-k,-m}^\ast(-\tau)\right).
\end{equation*}
\end{theorem}

To prove Theorem~\ref{thm:integer-weight}, we determine the extension of the holomorphic and non-holomorphic parts of $F_{k,-m}$ separately, in Sections~\ref{sec:hol} and \ref{sec:nonhol}, respectively.
The computation involving the holomorphic part closely follows \cite{Rh}, and the extension is provided by the simple fact that
\[
	\frac{1}{1-q} = 
	\begin{dcases}
		\sum_{n\geq 0} q^n & \text{ if }|q|<1, \\
		-\sum_{n \geq 1} q^{-n} & \text{ if }|q|>1.
	\end{dcases}
\]
For the non-holomorphic part, the situation is similar, but somewhat more complicated, and the extension is provided by the functional equation of the polylogarithm $\Li_{1-k}(q)$ (defined in \eqref{eq:Li-def} below), namely

\begin{equation} \label{eq:Li-fun-intro}
	\Li_{k-1} (q) = \Li_{k-1}\left(q^{-1}\right) \qquad \text{ for }   k\in -2\N.
\end{equation}

\begin{remark}
If one tries to mimic the computations of Section~\ref{sec:nonhol} in the case of half-integral weight, the situation is complicated by the analogue of \eqref{eq:Li-fun-intro} for $k\notin\Z$, namely
\begin{equation*} 
	i^{1-k} \Li_{k-1}\left(e^{2\pi i u}\right) + i^{k-1} \Li_{k-1}\left(e^{-2\pi i u}\right) = \frac{(2\pi)^{k-1}}{\Gamma(k-1)} \zeta(2-k,u),
\end{equation*}
where $\zeta(2-k,u)$ denotes the Hurwitz zeta function.
It is unclear whether the resulting function in the lower-half plane has any relation to a known modular-type object.
\end{remark}


The paper is organized as follows.
In Section~\ref{sec:preliminaries} we recall the definitions and some basic properties of harmonic Maass forms and Poincar\'e series.
In Section~\ref{sec:proof} we prove Theorem~\ref{thm:integer-weight}.


\section{Preliminaries} \label{sec:preliminaries}

\subsection{Harmonic Maass forms}
In this section we recall basic facts of harmonic Maass forms, first introduced by Bruinier and Funke in \cite{BF}. We begin with their definition.
\begin{definition}
For $k\in2\N$, a {\it weight $k$ harmonic Maass form for $\Gamma_0(N)$} is any smooth function $f:\mathbb H\rightarrow \C$ satisfying the following conditions:
\begin{enumerate}
	\item For all $\left(\begin{smallmatrix}
		a & b \\ c & d
	\end{smallmatrix} \right)\in \Gamma_0(N)$ we have
	$$
	f\left(\frac{a\tau+b}{c\tau+d}\right)=(c\tau+d)^kf(\tau).
	$$
	\item We have $\Delta_k(f)=0$ where $\Delta_k$ is the {\it weight $k$ hyperbolic Laplacian}
	$$
	\Delta_k:=-v^2\left(\frac{\partial^2}{\partial u^2}+\frac{\partial^2}{\partial v^2}\right)+ikv\left(\frac{\partial}{\partial u}+i\frac{\partial}{\partial v}\right).
	$$
	\item There exists a polynomial $P_f(\tau)\in\C [q^{-1}]$ such that
	$$
	f(\tau)-P_f(\tau)=O(e^{\varepsilon v})
	$$
	as $v\to \infty$ for some $\varepsilon >0$. Analogous conditions are required at all cusps.
\end{enumerate}
\end{definition}
Denote the space of such harmonic Maass forms by $H_k\left(\Gamma_0(N)\right)$. Every $f\in H_k\left(\Gamma_0(N)\right)$ has an expansion of the form
$$
f(\tau)=f^+(\tau)+f^-(\tau)
$$
with the {\it holomorphic part} having a $q$-expansion
\begin{equation*}
f^+(\tau)=\sum_{n\gg -\infty}c^+_f(n)q^n
\end{equation*}
and the {\it non-holomorphic part} having an expansion of the form
\begin{equation*}
f^-(\tau)=\sum_{n> 0}c_f^-(n)\Gamma(1-k,4\pi nv)q^{-n}.
\end{equation*}
Here $\Gamma(s,v)$ is the {\it incomplete gamma function} defined, for $v>0$, as the integral
\begin{equation}\label{incomplete}
	\Gamma(s,v) := \int_{v}^\infty t^{s-1}e^{-t} \, dt.
\end{equation}
%

\subsection{Poincar\'e series}

In this section, we recall the definitions and properties of various Poincar\'e series. The general construction is as follows. Let $\varphi$ be any translation-invariant function, which we call the \emph{seed} of the Poincar\'e series in question. Then, in the case of absolute convergence, we can define a function satisfying weight $k$ modularity by forming the sum 
\begin{equation*}
\mathbb{P}_{k}(\varphi;\tau):=\sum_{\gamma \in  \Gamma_{\infty}  \backslash \Gamma_0(N)} \varphi|_k\gamma (\tau),
\end{equation*}
where $\Gamma_{\infty}:= \{\pm \left(\begin{smallmatrix}1&n\\0&1 \end{smallmatrix} \right)   : n \in \mathbb{Z}  \}$ is the group of translations. Convergence is, in particular, satisfied by functions $\varphi$ satisfying
$\varphi(\tau)= O (v^{2-k+\varepsilon})$ as $v\to 0$.

A natural choice for $\varphi$ is a typical Fourier coefficient in the space of automorphic functions one is interested in. For example,
in the case of weakly holomorphic modular forms one may choose, for $m \in \mathbb{Z}$,
$$
\varphi(\tau)=\varphi_m(\tau):= 	q^m.
$$
Define for $k\in2\N$ with $k > 2$ and $m\in\mathbb{Z}$ the {\it Poincar\'e series of exponential type} by
\begin{equation*}  
	P_{k,m}(\tau)
	:=\mathbb{P}_{k}(\varphi_m;\tau)=\sum_{\gamma \in  \Gamma_{\infty}  \backslash \Gamma_0(N)} \varphi_m|_k\gamma (\tau).
\end{equation*}
To give their Fourier expansion, we require the \emph{Kloosterman sums}
\begin{equation} 
\label{Klooster_ind}
K(m,n;c):=
\sum_{d\pmod c^*}  e\left(\frac{m\overline d+nd}{c}\right),
\end{equation}
where $e(x) := e^{2\pi i x}$.
A direct calculation yields the following duality:
\begin{equation*}
	K(-m,-n;c) = K(m,n;c).
\end{equation*}

A very useful property of the Poincar\'e series is that they have explicit Fourier expansions, as given in the following theorem.
\begin{theorem}\label{FourierExpa}
Suppose that $k>2$ is even.
	\begin{itemize}
		\item[i)]
		If $m\in\N$, the Poincar\'e series $P_{k,m}$ are in $S_k(\Gamma_0(N))$.
		 We have the Fourier expansion $P_{k,m} (\tau)  =\sum_{n=1}^\infty b_{k,m}(n)q^n,$ where
		$$
		b_{k,m}(n)
		= \left(\frac{n}{m} \right)^{\frac{k-1}{2}} \left( \delta_{m,n} + 2 \pi (-1)^{\frac k2} \sum_{ \substack{c>0 \\ N\mid c}} \frac{K(m,n;c)}{c} J_{k-1} \left( \frac{4 \pi \sqrt{mn}}{c}\right)\right).
		$$
        Here $\delta_{m,n}$ is the Kronecker delta-function and $J_s$ denotes the usual $J$-Bessel function.
		\item[ii)]
		For $m\in-\mathbb{N}$, the  Poincar\'e series $ P_{k,m} $ are elements of  $M_k^!(\Gamma_0(N))$.  We have the Fourier expansion $ {P}_{k,m} (\tau)  =q^{m} +\sum_{n=1}^{\infty} b_{k,m}(n)q^n,$ where
		$$
		b_{k,m}(n)
		= 2 \pi (-1)^{\frac k2}\left|\frac{n}{m} \right|^{\frac{k-1}{2}}
		\sum_{ \substack{c>0 \\ N\mid c}}
		\frac{K(m,n;c)}{c} I_{k-1}
		\left(
		\frac{4 \pi \sqrt{|mn|}}{c}
		\right).
		$$
		Here $I_s$ denotes the usual I-Bessel function.
		Moreover, $P_{k,m}$ is holomorphic at the cusps of $\Gamma_0(N)$ other than $i\infty$.
	\end{itemize}
\end{theorem}	
We next turn to the construction of harmonic Maass forms via Poincar\'e series.   Such series have appeared in many places
in the literature, indeed in the works of Niebur \cite{niebur} and Fay \cite{fay} in the 1970's, long before the recent advent of harmonic Maass forms.
Define 
	\begin{equation*}
	F_{k,m}
	:=\sum_{\gamma\in  \Gamma_{\infty}  \backslash \Gamma_0(N)} \phi_{k,m}|_k\gamma,
	\end{equation*}
where the seed $\phi_{k,m}$ is given by
\[
	\phi_{k,m}(\tau) :=  \big(1 - \Gamma^*(1-k,4\pi |m| v) \big) q^m.
\]
Here $\Gamma^*$ is the {\it normalized incomplete gamma function}
\begin{equation*}
\Gamma^*(s,v) := \frac{\Gamma(s,v)}{\Gamma(s)}.
\end{equation*}
The analogous exact formula for coefficients of these Poincar\'e series is then given in the following theorem (see, e.g., \cite{fay} or \cite{hmf-book} for a proof).
\begin{theorem}\label{MaassPoincareThm}
	If $k < 0$ is even and $m\in-\N$, then $F_{k,m} \in H_{k}(\Gamma_0(N))$.
	We have
\begin{equation*}
	\xi_k\left(F_{k,m}\right) = -\frac{(4\pi m)^{1-k}}{(-k)!} P_{2-k,-m}
\end{equation*}
and
\begin{equation*}
	D^{1-k}\left(F_{k,m}\right) = m^{1-k}P_{2-k,m}.
\end{equation*}
We have the Fourier expansion
\begin{multline*}
	F_{k,m}(\tau) =\left(1-\Gamma^*(1-k,-4\pi mv)\right) q^m \\
	+ \sum_{n = 0}^\infty a_{k,m}^+(n)q^n + \sum_{n =1}^\infty a_{k,m}^-(n)\Gamma^*(1-k,4\pi n v)q^{-n}
\end{multline*}
with
\[
	a_{k,m}^+(0) = \frac{(2\pi)^{2-k}(-1)^{\frac k2+1}m^{1-k}}{(1-k)!}\sum_{\substack{c > 0 \\ N\mid c}} \frac{K(m,0;c)}{c^{2-k}}.
\]
Moreover, for $ n \geq 1$ and $\varepsilon \in \{+, -\}$, we have 

\begin{multline*}
a_{k,m}^\varepsilon (n)
= 2\pi (-1)^{\frac{k}{2}}  \left\lvert\frac{m}{n}\right\rvert^{\frac{1-k}{2}}
\sum_{\substack{c>0 \\ N\mid c}}\frac{K (m, \varepsilon n; c)}{c} \times
\begin{cases}
I_{1-k}\left(\frac{4\pi\sqrt{|mn|}}{c}\right) &\text{ if } \varepsilon n >0,\\\vspace*{-4mm}\\
J_{1-k}\left(\frac{4\pi\sqrt{|mn|}}{c}\right) &\text{ if } \varepsilon n <0.
\end{cases}
\end{multline*}
\end{theorem}

\section{Proof of Theorem \ref{thm:integer-weight}}
\label{sec:proof}

To prove Theorem \ref{thm:integer-weight}, we consider the holomorphic and non-holomorphic parts of $F_{k,-m}$ separately.

\subsection{The holomorphic part}
\label{sec:hol}

We first extend the holomorphic part $F_{k,-m}^+$ of $F_{k,-m}$
to a function defined for $|q|\neq 1$, closely following \cite{Rh}. Using Theorem~\ref{MaassPoincareThm},
we have, for $|q|<1$,
\begin{multline*}
	F_{k,-m}^+(\tau) \\
	= q^{-m} + a_{k,-m}^+(0) + 2\pi (-1)^{\frac k2} m^{\frac{1-k}2} \sum_{\substack{c>0 \\ N\mid c}} \frac 1c \sum_{d\pmod{c}^*} e\pfrac{-m\bar d}{c} A^+_m(c,d),
\end{multline*}
where
\begin{equation} \label{eq:A+-def}
	A^+_m(c,d) = A^+_m(c,d;\tau) := \sum_{n\geq 1} n^{\frac{k-1}2}I_{1-k}\pfrac{4\pi \sqrt{mn}}{c} \zeta_c^{nd} q^n.
\end{equation}
Using the series expansion of the $I$-Bessel function
\[
	I_\alpha(x)=\sum_{j\geq 1}\frac{1}{j!\Gamma(j+\alpha+1)}\left(\frac{x}{2}\right)^{2j+\alpha},
\]
we obtain
$$
	n^{\frac{k-1}2} I_{1-k} \pfrac{4\pi \sqrt{mn}}{c}= \sum_{j\geq 0} \beta^+_{m,c}(j) \frac{n^j}{j!},
$$
where
\begin{equation*}
	\beta^+_{m,c}(j) := \frac{\pfrac{2\pi\sqrt{m}}{c}^{2j+1-k}}{(j+1-k)!}.
\end{equation*}
We insert the integral representation (for $r>0$)
\begin{equation}
\label{rewritefak}
	\frac{n^j}{j!} = \frac{1}{2\pi i} \int_{|s|=r} \frac{e^{ns}}{s^{j+1}} \, ds,
\end{equation}
and we conclude that
\begin{equation*} 
	n^{\frac{k-1}2} I_{1-k} \pfrac{4\pi \sqrt{mn}}{c}
	= \frac{1}{2\pi i}\int_{|s|=r} \alpha^+_{m,c}(s) e^{ns} \, ds,
\end{equation*}
where $\alpha^+_{m,c}(s)$ is the series
\[
	\alpha^+_{m,c}(s) := \sum_{j\geq0} \frac{\beta^+_{m,c}(j)}{s^{j+1}},
\]
which is absolutely convergent for all $s$.
Equation \eqref{eq:A+-def} then becomes
\begin{align*}
	A^+_m(c,d)
	= \frac{1}{2\pi i} \int_{|s|=r} \alpha^+_{m,c}(s) \sum_{n\geq 1} \left(e^s \zeta_c^d q\right)^n \, ds
	= \frac{1}{2\pi i} \zeta_c^d q\int_{|s|=r} \frac{\alpha^+_{m,c}(s)e^s}{1-e^s \zeta_c^d q} \, ds.
\end{align*}
Here we take $r$ sufficiently small so that $|e^s \zeta_c^d q|<1$.
Define
\[
	\phi_{k,m}^+(c,d;\tau) := \frac{1}{2\pi i} \int_{|s|=r} \frac{\alpha^+_{m,c}(s)e^s}{1-e^s\zeta_c^d q} \, ds.
\]
We can now define the function which exists away from the real line.
To be more precise, since $\phi^+_{k,m}$ is regular for all $v\neq 0$, the function
\begin{multline} \label{eq:H-k-m-+}
	H_{k,m}^+(\tau) := q^{-m} + a_{k,-m}^+(0) \\
	+ 2\pi (-1)^{\frac k2} m^{\frac{1-k}2} \sum_{\substack{c>0 \\ N\mid c}} \frac 1c \sum_{d\pmod{c}^*} e\left(\frac{-m\bar d+d}{c}+\tau\right) \phi_{k,m}^+(c,d;\tau)
\end{multline}
is defined for $\tau\in \H \cup (-\H)$.

We now consider the Fourier expansion of the function $H_{k,m}^+(\tau)$ for $\tau$ in the lower half-plane, so suppose for the remainder of the proof that $v<0$. In this case, we have
\begin{align*}
 	\zeta_c^d q \, \phi_{k,m}^{+}(c,d;\tau)
 	&= \frac{1}{2\pi i} \int_{|s|=r} \alpha^+_{m,c}(s) \frac{e^s \zeta_c^d q}{1-e^s \zeta_c^d q} ds \\
 	&= -\frac{1}{2\pi i} \int_{|s|=r} \alpha^+_{m,c}(s) \sum_{n\geq 0} \left(e^s \zeta_c^d q\right)^{-n} ds,
\end{align*}
where $r$ is chosen so that $|e^{-s}q^{-1}|<1$.
By reversing the calculation which led to \eqref{eq:H-k-m-+}, making the change of variables $s \mapsto -s$, and using  that $I_{1-k}(-ix)= i^{1-k}J_{1-k}(x)$, we find that
\begin{align*}
 	\zeta_c^d q \, \phi_{k,m}^{+}(c,d;\tau) &= -\beta^+_{m,c}(0) -\sum_{n\geq 1} \left(\sum_{j\geq 0} (-1)^j \beta^+_{m,c}(j) \frac{n^j}{j!}\right) e\pfrac{-nd}{c} q^{-n} \\
 	&= -\beta^+_{m,c}(0) -\sum_{n\geq 1} e\pfrac{-nd}{c} n^{\frac{k-1}2} J_{1-k}\pfrac{4\pi\sqrt{mn}}{c} q^{-n}.
\end{align*}
Thus we have
\begin{multline} \label{eq:H-plus-outside-1}
 	H_{k,m}^+(\tau) \\
 	=q^{-m} + a_{k,-m}^+(0) - 2\pi (-1)^{\frac k2} m^{\frac{1-k}2}  \sum_{\substack{c>0 \\ c\equiv 0 \pmod{N}}} \frac 1c \sum_{d \pmod{c}^*} e\pfrac{-m\bar d}{c} \beta^+_{m,c}(0) \\
 	- 2\pi (-1)^{\frac k2} \sum_{n\geq 1}\left( \pfrac{n}{m}^{\frac{k-1}{2}} \sum_{\substack{c>0\\N\mid c}} \frac{K(-m,-n;c)}{c} J_{1-k}\pfrac{4\pi\sqrt{mn}}{c} \right) q^{-n}.
\end{multline}
Note that the second and the third terms  on the right-hand side of \eqref{eq:H-plus-outside-1} cancel, since
\begin{align*}
	- 2\pi (-1)^{\frac k2} & m^{\frac{1-k}2}  \sum_{\substack{c>0 \\ N\mid c}} \frac 1c \sum_{d\pmod{c}^*}  e\pfrac{-m\bar d}{c} \beta^+_{m,c}(0)\\
	&= \frac{(2\pi)^{2-k}(-1)^{\frac k2+1}m^{1-k}}{(1-k)!} \sum_{\substack{c>0 \\ N\mid c}} \frac{K(-m,0;c)}{c^{2-k}} = -a_{k,-m}^+(0).
\end{align*}
By (\ref{Klooster_ind}), we conclude that
\begin{align*}
	H_{k,m}^+(\tau)
	&= q^{-m} - \sum_{n\geq 1} \pfrac{m}{n}^{1-k}\\
	&\quad\times\left( 2\pi (-1)^{\frac k2} \pfrac{n}{m}^{\frac{1-k}{2}} \sum_{\substack{c>0\\N\mid c}} \frac{K(m,n;c)}{c} J_{1-k}\pfrac{4\pi\sqrt{mn}}{c} \right) q^{-n} \displaybreak[0] \\
	&= m^{1-k}\left( m^{k-1}q^{-m} + \sum_{n\geq 1} n^{k-1} b_{2-k,m}(n) q^{-n} \right)= m^{1-k} \mathcal E_{P_{2-k,m}}(-\tau)
\end{align*}
if $\tau$ is in the lower half-plane.

\subsection{The non-holomorphic part}
\label{sec:nonhol}
Next we extend the non-holomorphic part $F_{k,-m}^-(\tau)$
to a function $H_{k,m}^-(\tau)$, which is defined for $|q|\neq 1$.
We have, by Theorem \ref{MaassPoincareThm},
\begin{multline*}
	F_{k,-m}^-(\tau)
	=-\Gamma^*(1-k,4\pi m v)q^{-m} \\ + 2\pi (-1)^{\frac k2}m^{\frac{1-k}2}\sum_{\substack{c>0 \\ N\mid c}} \frac 1c \sum_{d\pmod{c}^*} e\left(\frac{-m\bar d}{c}\right) A_m^-(c,d;\tau),
\end{multline*}
where
\begin{align*}
	A_m^-(c,d) &= A_m^-(c,d;\tau) \\
	&:= \sum_{n\geq 1} n^{\frac{k-1}2} e\pfrac{-nd}{c} J_{1-k}\pfrac{4\pi\sqrt{mn}}{c} \Gamma^*(1-k,4\pi nv) q^{-n}.
\end{align*}
Using the integral representation of the incomplete Gamma function and making a change of variables, we find that
\[
	\Gamma(1-k,y) = y^{1-k}\int_1^\infty t^{-k} e^{-yt} dt,
\]
thus
\begin{align*}
	A_m^-(c,d) = \frac{(4\pi v)^{1-k}}{(-k)!} \int_1^\infty t^{-k} \sum_{n\geq 1} n^{\frac{1-k}2} J_{1-k} \pfrac{4\pi \sqrt{mn}}{c} e^{-4\pi n v t} \left(\zeta_c^d q\right)^{-n} dt.
\end{align*}
As above, we use the series expansion of the $J$-Bessel function 
\[
	J_\alpha(x)=\sum_{j\geq 1}\frac{(-1)^j}{j!\Gamma(j+\alpha+1)}\left(\frac{x}{2}\right)^{2j+\alpha}
\]
to obtain
\begin{align*}
	n^{\frac{1-k}2} J_{1-k} \pfrac{4\pi\sqrt{mn}}{c}
	&= n^{1-k} \sum_{j\geq 0} \beta_{m,c}^{-}(j) \frac{n^j}{j!},
\end{align*}
where
\[
	\beta_{m,c}^-(j) := \frac{(-1)^j}{(j+1-k)!}\pfrac{2\pi\sqrt{m}}{c}^{2j+1-k}.
\]
Thus we have, again using \eqref{rewritefak},
\begin{align*}
	n^{\frac{1-k}2} J_{1-k} \pfrac{4\pi\sqrt{mn}}{c}
	&= \frac{1}{2\pi i} \int_{|s|=r} n^{1-k} e^{ns} \alpha_{m,c}^{-}(s) ds,
\end{align*}
where
\[
	\alpha_{m,c}^-(s) := \sum_{j\geq 0} \frac{\beta_{m,c}^-(j)}{s^{j+1}}.
\]
Here $r$ is chosen so that $r<2\pi v$.
Thus
\begin{align*}
	A_m^-(c,d) & = \frac{(4\pi v)^{1-k}}{(-k)!} \int_1^\infty t^{-k} \left( \frac{1}{2\pi i} \int_{|s|=r} \alpha_{m,c}^-(s) \sum_{n\geq 1} n^{1-k} \left(e^{s-4\pi  v t} \zeta_c^{-d} q^{-1}\right)^n ds \right) dt\\
                   & = \frac{(4\pi v)^{1-k}}{(-k)!} \int_1^\infty t^{-k} \left( \frac{1}{2\pi i} \int_{|s|=r} \alpha_{m,c}^-(s) \Li_{k-1}\left(e^{s-4\pi vt}\zeta_c^{-d} q^{-1}\right) ds \right) dt.
\end{align*}
Here $\Li_s(w)$ is the {\it polylogarithm}, defined for $s\in \C$ and $|w|<1$ by
\begin{equation} \label{eq:Li-def}
	\Li_s(w) := \sum_{n=1}^\infty \frac{w^n}{n^s}.
\end{equation}

We now again introduce a function defined away from the real line. Define
\begin{equation*}
	\phi_{k,m}^-(c,d;\tau) := \frac{1}{2\pi i} \int_1^\infty t^{-k} \int_{|s|=r} \alpha_{m,c}^-(s) \Li_{k-1}\left(e^{s-4\pi vt}\zeta_c^{-d}q^{-1}\right) \, ds \, dt
\end{equation*}
and
\begin{multline} \label{eq:h-k-m--}
	H_{k,m}^-(\tau) := -\Gamma^*(1-k,4\pi mv)q^{-m} \\ + \frac{2\pi(-1)^{\frac k2}(4\pi v)^{1-k}}{(-k)!} m^{\frac{1-k}2} \sum_{\substack{c>0\\N\mid c}} \frac 1c \sum_{d\pmod{c}^*} e\pfrac{-m\bar d}{c} \phi_{k,m}^-(c,d;\tau).
\end{multline}
Using the functional equation 
\begin{equation*}
\mathrm{Li}_{-n}(z)=(-1)^{n+1}\mathrm{Li}_{-n}\left(\frac{1}{z}\right),
\end{equation*}
we obtain, for $v<0$, and using that $k$ is even,
\begin{equation*}
	A_m^-(c,d) = \frac{(4\pi v)^{1-k}}{(-k)!} \int_1^\infty t^{-k} \left( \frac{1}{2\pi i} \int_{|s|=r} \alpha_{m,c}^-(s) \Li_{k-1}\left(e^{-s+4\pi vt}\zeta_c^d q\right) \, ds \right)  dt.
\end{equation*}
Since $v<0$ we can now use the series representation of $\Li_{k-1}$.
This yields
\begin{equation*}
	A_m^-(c,d) = \frac{(4\pi v)^{1-k}}{(-k)!} \! \int_1^\infty t^{-k} \left( \sum_{n\geq 1} n^{1-k} e^{4\pi n v t} \left(\zeta_c^d q\right)^n \frac{1}{2\pi i} \int_{|s|=r} \! \alpha_{m,c}^-(s) e^{-ns} \, ds \! \right) \! dt.
\end{equation*}
The innermost integral is (inserting the definition of  $\alpha_{m,c}^-(s)$, making the change of variables $s \mapsto -s$, and using \eqref{rewritefak})
\begin{align*}
	\frac{1}{2\pi i} \int_{|s|=r} \alpha_{m,c}^-(s) e^{-ns} ds
	&= \sum_{j\geq 0} \frac{1}{(j+1-k)!} \pfrac{2\pi \sqrt{m}}{c}^{2j+1-k} \frac{n^j}{j!}\\
	&= n^{\frac{k-1}2} I_{1-k} \pfrac{4\pi\sqrt{mn}}{c}.
\end{align*}
Thus
\begin{align*}
	A_m^-(c,d) 	&= \frac{ (4\pi v)^{1-k}}{(-k)!} \sum_{n\geq 1} n^{\frac{1-k}2} I_{1-k}\pfrac{4\pi\sqrt{mn}}{c} \left( \int_1^\infty t^{-k}e^{4\pi n v t} dt \right) \left(\zeta_c^d q\right)^n \\
	&= \sum_{n\geq 1} n^{\frac{k-1}2} I_{1-k}\pfrac{4\pi\sqrt{mn}}{c} \Gamma^*(1-k,4\pi n |v|) \left(\zeta_c^d q\right)^n.
\end{align*}
Therefore, using that $K(m, n; c)$ is real, we conclude that
\begin{align*}
    H_{k,m}^-(\tau) = &-\Gamma^*(1-k,4\pi mv)q^{-m} 
    +2\pi (-1)^{\frac k2} \sum_{n\geq 1} \pfrac{n}{m}^{\frac{k-1}2} \\
    &\qquad \quad \times \sum_{\substack{c>0 \\ N\mid c}} \frac{K(-m,n;c)}{c} I_{1-k}\pfrac{4\pi\sqrt{mn}}{c} \Gamma^*(1-k,4\pi n|v|) q^n \\
	&= -\frac{(4\pi m)^{1-k}}{(-k)!} P^{*}_{2-k,-m}(-\tau).
\end{align*}


\bibliographystyle{plain}

\end{document}